\documentclass[a4paper,12pt,reqno]{amsart}
\usepackage{amssymb,delarray,amsmath,stmaryrd}
\usepackage[dvips]{graphicx}

\let\3=\ss

\let\ss=\sss

\usepackage{hyperref}
\usepackage[all]{xy}

\input epsf.tex
\newdimen\xsize
\newdimen\oldbaselineskip
\newdimen\oldlineskiplimit
\def\restorelineskip{\baselineskip=\oldbaselineskip%
\lineskiplimit=\oldlineskiplimit}
\def\putm[#1][#2]#3{
\hbox{\vbox to 0pt{\parindent=0pt%
\vskip#2\xsize\hbox to0pt{\hskip#1\xsize $#3$\hss}\vss}}}%
\long\def\Line#1{\hbox to \hsize{#1}}
\def\putt[#1][#2]#3{
\vbox to 0pt{\noindent\hskip#1\xsize\lower#2\xsize%
\vtop{\restorelineskip#3}\vss}}

\makeatletter
\def\xbig[#1]#2{{\hbox{$\m@th\left#2\vbox to#1\xsize{}%
\right.\n@space$}}}
\makeatother
\def\xlar[#1]#2{%
\smash{\mathop{ \hbox to #1\xsize{\leftarrowfill}}\limits^{#2}}}
\def\xrar[#1]#2{%
\smash{\mathop{ \hbox to #1\xsize{\rightarrowfill}}\limits^{#2}}}
\def\xline[#1]{\hbox to #1\xsize{\leaders\hrule\hfill}}

\DeclareFontFamily{U}{rsf}{\skewchar\font'177}%
\DeclareFontShape{U}{rsf}{m}{n}{<-6>rsfs5<6-8>rsfs7<8->rsfs10}{}%
\DeclareFontShape{U}{rsf}{b}{n}{<-6>rsfs5<6-8>rsfs7<8->rsfs10}{}%
\DeclareMathAlphabet\RSFS{U}{rsf}{m}{n}
\SetMathAlphabet\RSFS{bold}{U}{rsf}{b}{n}
\def\sf#1{{\mathsf{#1}}}

\def\slsf{\slshape \sffamily }

\def\msmall#1{\mathchoice{\hbox{\small$\displaystyle {#1}$}}{#1}{#1}{#1}}

\def\ss{{\mathbb S}}

\def\cc{{\mathbb C}}

\def\rr{{\mathbb R}}

\def\im{\sf{Im}\,}

\def\re{\sf{Re}\,}

\def\lim{\mathop{\sf{lim}}}

\def\eps{\varepsilon}

\def\<{\langle}\let\la=\<
\def\>{\rangle}\let\ra=\>

\def\d{\partial}
\def\dbar{{\barr\partial}}

\def\ddef{\mathrel{{=}\raise0.3pt\hbox{:}}}
\def\deff{\mathrel{\raise0.3pt\hbox{\rm:}{=}}}

\def\fraction#1/#2{\mathchoice{{\msmall{ #1\over#2}}}%
{{ #1\over #2 }}{{#1/#2}}{{#1/#2}}}
\def\norm#1{\left\Vert{#1}\right\Vert}

\def\emptyset{\varnothing}

\def\longpoints{\leaders\hbox to 0.5em{\hss.\hss}\hfill \hskip0pt}
\def\stateskip{\smallskip}
\def\state#1. {\stateskip\noindent{\bf#1. }} 
\def\statep#1. {\stateskip\noindent{\bf#1 }} 
\def\proof{\state Proof. \,}

\def\Chi{\raise 2pt\hbox{$\chi$}}

\def\ie{\hskip1pt plus1pt{\sl i.e.\/,\ \hskip1pt plus1pt}}

\def\sli{{\sl i)} } 
\def\slii{{\sl i$\!$i)} } 
\def\sliii{{\sl i$\!$i$\!$i)} }

\def\barr#1{\mskip1mu\overline{\mskip-1mu{#1}\mskip-1mu}\mskip1mu}
\def\Chi{\raise 2pt\hbox{$\chi$}}

\let\phI=\phi\let\phi=\varphi\let\varphi=\phI
\let\cal=\mathcal

\def\calc{{\cal C}}

\def\calp{{\cal P}}

\def\eps{\varepsilon}

\def\d{\partial}
\def\dbar{{\barr\partial}}
\def\1{{1\mkern-5mu{\rom l}}}

\def\ge{\geqslant}

\def\fraction#1/#2{\mathchoice{{\msmall{ #1\over#2}}}%
{{ #1\over #2 }}{{#1/#2}}{{#1/#2}}}

\def\emptyset{\varnothing}

\def\qed{\ \ \hfill\hbox to .1pt{}\hfill\hbox to .1pt{}\hfill $\square$\par}

\def\comment#1\endcomment{}

 \belowdisplayskip=\abovedisplayskip

\def\lineeqqno(#1){\hfill\llap{\vbox to 10pt%
{\vss\begin{align} \eqqno(#1)\end{align}\vss}}\vskip1pt}

\advance\headheight 1.2pt

\def\ShowwLLabel#1{}

\def\thechpt{\Roman{chpt}}

\def\newchapt[#1]#2{\newpage%
\refstepcounter{chpt}\setcounter{subsection}{0}%
\setcounter{thm}{0}\setcounter{defi}{0}%
\setcounter{rema}{0}\setcounter{exrc}{0}%
\renewcommand{\thesubsection}{\thechpt.\arabic{subsection}}%
\section*{\begin{center}\huge \bf Chapter \thechpt\\
#2 \end{center}}\label{#1}%
\ \smallskip%
\markboth{Chapter \thechpt}{#2}%
}
\def\newsect[#1]#2{\refstepcounter{section}\setcounter{equation}{0}%
\renewcommand{\thesubsection}{\arabic{section}.\arabic{subsection}}%
\section*{\arabic{section}.
#2}\vspace{-20pt}\label{#1}\vspace{20pt}%
\markboth{Section \arabic{section}}{#2}}

\def\newlect[#1]#2{\refstepcounter{section}%
\renewcommand{\thesubsection}{\arabic{section}.\arabic{subsection}}%
\section*{Lecture \arabic{section}\\
#2}\label{#1}%
\markboth{Lecture \arabic{section}}{#2}}

\def\newprg[#1]#2{\refstepcounter{subsection}%
\subsection*{{\thesubsection.\ #2}} \label{#1}%
}

\def\newappx[#1]#2{%
\refstepcounter{appx}\setcounter{section}{0}%
\renewcommand{\thesubsection}{A\arabic{appx}.\arabic{subsection}}%
\section*{Appendix \arabic{appx}\\ #2}
\label{#1}%
\markboth{Appendix A\arabic{appx}}{#2}
}

\newtheorem{thm}{Theorem}
   \def\newthm#1{\begin{thm}\label{#1}}

\newtheorem{nnthm}{Theorem.} 
   \def\newnnthm#1{\begin{nnthm} \label{#1}}
\newtheorem{lem}{Lemma}[section]
   \def\newlemma#1{\begin{lem} \label{#1}}

\newtheorem{prop}{Proposition}[section]
   \def\newprop#1{\begin{prop}\label{#1}}

\newtheorem{corol}{Corollary}[section]
   \def\newcorol#1{\begin{corol} \label{#1}}

\newtheorem{defi}{Definition}[section]
   \def\newdefi#1{\begin{defi} \label{#1}\rm }

\newtheorem{exmp}{Example}[section]
   \def\newexmp#1{\begin{exmp} \label{#1}\rm }

\newtheorem{exrc}{Exercise}
   \def\newexrc#1{\begin{exrc} \label{#1}\rm }

\newtheorem{quest}{Open Question}
   \def\newrema#1{\begin{quest} \label{#1}\rm }

\newtheorem{rema}{Remark}
   \def\newrema#1{\begin{rema} \label{#1}\rm }

\let\xrar=\xrightarrow

\def\eqqno(#1){\label{(#1)}}
\def\eqqref(#1){(\ref{(#1)})}

\numberwithin{equation}{section}

\title{On convex to pseudoconvex mappings}
\author{S. Ivashkovich}
\address{
Universit\'e de Lille-1, UFR de Math\'ematiques, 59655 Villeneuve
d'Ascq, France.}
\email{ivachkov@math.univ-lille1.fr}
\address{IAPMM Nat. Acad. Sci. Ukraine,
Lviv, Naukova 3b,
79601 Ukraine.}
\date{\today}
\begin{document}
\maketitle

\begin{abstract}
In the works of Darboux and Walsh, see \cite{D,W}, it was remarked that a one to one self mapping of $\rr^3$ which 
sends convex sets to convex ones is affine. It can be remarked also that a $\calc^2$-diffeomorphism $F:U\to U^{'}$ between 
two domains in $\cc^n$, $n\ge 2$, which sends pseudoconvex hypersurfaces to pseudoconvex ones  is either holomorphic or 
antiholomorphic. 

\smallskip In this note we are interested in the self mappings of $\cc^n$ which send convex hypersurfaces to pseudoconvex
ones. Their characterization is the following: {\it A $\calc^2$ - diffeomorphism $F:U'\to U$ (where $U', U\subset \cc^n$ are
domains) sends convex hypersurfaces to pseudoconvex ones if and only if the inverse map $\Phi\deff F^{-1}$ is 
weakly pluriharmonic, \ie it satisfies some nice second order PDE very close to $\d\bar\d \Phi = 0$.} In fact all pluriharmonic
$\Phi$-s do satisfy this equation, but there are also other solutions.
\end{abstract}

\section[1]{Formulation} Let $U', U$ be domains in $\cc^n, n\ge 2$ and let $F: U'\to U$ be a $\calc^2$-diffeomorphism. 
Coordinates in the source we denote by $z'=x'+iy'$, in the target by $z=x+iy$. It will be convenient for us to 
suppose that $U'$ is a convex neighborhood of zero and that $F(0')=0$. The, somewhat unusual choice to put primes
on the objects in the source (and not in the target) is explained by the fact that in the statements and in the proofs 
we shall work more with the inverse map $\Phi$ then with $F$.

\begin{thm}
Let $F:U'\to U$ be a $\calc^2$-diffeomorphism. Then the following conditions are equivalent:

\sli For every convex hypersurface $M'\subset U'$ the image $M=F(M')$ is a
pseudoconvex hypersurface in $U$.

\slii The inverse map $\Phi\deff F^{-1}:U\to U'$ satisfies the following second order PDE System
\begin{equation}
\eqqno(syst1)
\d\dbar \Phi = (d\Phi^{-1}(\Delta\Phi),dz)\wedge \d\Phi + (dz,d\Phi^{-1}(\Delta \Phi))\wedge\dbar \Phi .
\end{equation}

\sliii The equation \eqqref(syst1) has the following geometric meaning: for every $z\in U$ and every $\zeta\in T_z\cc^n$
\begin{equation}
\eqqno(syst2)
\d\dbar\Phi_z (\zeta ,\bar\zeta) \in \texttt{span} \left\lbrace d\Phi_z(\zeta),d\Phi_z(i\zeta)\right\rbrace .
\end{equation}
\end{thm}
Here we use the following notation: for a vector $v=(v^1,...,v^n)\in \cc^n$ and $dz=(dz_1,..., dz_n)$ we set 
$(dz,v) = \bar v^jdz_j$ and $(v,dz) = v^jd\bar z_j$. Throughout this note we shall use the Einstein summation convention. 

\begin{rema}\rm
Pluriharmonic $\Phi$-s clearly satisfy \eqqref(syst1) (or \eqqref(syst2)) and let us remark that this 
geometric characterization of pluriharmonic diffeomorphisms perfectly agrees with an analytic one:
{\it The class $\calp$ of  pluriharmonic diffeomorphisms $\cc^n \to \cc^n$ is stable under biholomorphic 
parametrization of the source and $\rr$-linear of the target}.
Really, these parametrization preserve accordingly pseudoconvexity and convexity of hypersurfaces.

\smallskip\noindent{\bf 2.} The item (\sli of the Theorem is clearly equivalent to the following one: {\it For 
every strictly convex quadric  $M'\cap U'\not= \emptyset$ the image $M=F(M'\cap U')$ is a pseudoconvex hypersurface 
in $U$}. I.e, it is enough to  check this condition only for  quadrics.

\smallskip\noindent{\bf 3.} The fact that \eqqref(syst1) admits other solutions then just pluriharmonic mappings
is very easy to see from the form of its linearization at identity
\begin{equation}
\eqqno(syst3)
\d\dbar\Phi = (\Delta\Phi , dz)\wedge dz .
\end{equation}
Remark that any map of the form $\Phi (z) = (\phi_1(z_1),...,\phi_n(z_n))$ satisfies \eqqref(syst3) provided all 
$\phi_j$, except for some $j_0$, are harmonic. And this $\phi_{j_0}$ can be then an arbitrary $\calc^2$-function.
\end{rema}

\section[2]{An auxiliary computation}

Denote by $\zeta = \xi + i\eta$ a tangent vector at point $z\in\cc^n$.  Recall that the  real Hessian of a real valued 
function $\rho$ in $\cc^n=\rr^{2n}$ is

\begin{equation}
\eqqno(r-hess)
H^{\rr}_{\rho (z)}(\zeta,\zeta) = \frac{\d^2\rho (z)}{\d x_i\d x_j}\xi_i\xi_j + \frac{\d^2\rho (z)}{\d y_i\d y_j}\eta_i\eta_j + 2\frac{\d^2\rho (z)}{\d x_i\d y_j}\xi_i\eta_j.
\end{equation}
A hypersurface  $M=\{ z\in U: \rho (z)=0\}$, with $\rho$ is $\calc^2$-regular, $\rho (0)=0$ and $\nabla\rho|_M\not= 0$, is strictly convex if the
defining function $\rho$ can be chosen with positive definite Hessian, \ie $H^{\rr}_{\rho (z)}(\zeta,\zeta)>0$
for all $z\in M$ and all $\zeta\not= 0$. One readily checks the following expression of the real Hessian of $\rho$ in complex 
coordinates

\begin{equation}
\eqqno(rc-hess)
H^{\rr}_{\rho (z)}(\zeta ,\zeta) = \frac{\d^2\rho (z)}{\d z_i\d z_j}\zeta_i\zeta_j +
\frac{\d^2\rho (z)}{\d\bar z_i\d\bar z_j}\bar\zeta_i\bar\zeta_j +
2\frac{\d^2\rho (z)}{\d z_i\d\bar z_j}\zeta_i\bar\zeta_j.
\end{equation}

Recall that the Hermitian part $L_{\rho (z)}(\zeta,\bar\zeta)=\frac{\d^2\rho}{\d z_i\d \bar z_j}\zeta_i\bar\zeta_j$ of the Hessian is called the Levi form of $\rho$ (and of $M$). $M$ is strictly
pseudoconvex if its Levi form is positive definite on the complex tangent space $T^c_zM=\{\zeta\in T_z\cc^n:
\left(\dbar\rho (z),\zeta \right)=0\}$ for every $z\in M$. Here $\left(\cdot,\cdot\right)$ stands for the standard
Hermitian scalar product in $\cc^n$. 


\smallskip Let $F: \cc^n_{z'}\supset U'\to U\subset \cc^n_{z}$ be a $\calc^2$-diffeomorphism. Let further $z'=z'(z)$ 
be the coordinate representation of the inverse mapping $z' = \Phi (z) \deff F^{-1} (z)$ and let
$M=F(M')\subset U$ be the image of a hypersurface $M'\subset U'$. Then $M =\{ z: \rho (z)=0\}$,
where $\rho (z):=\rho' (z'(z))$. 

\begin{lem}
The Levi form of $\rho $ at point $z$ decomposes as
\begin{equation}
\eqqno(levif-1)
L_{\rho (z)}(\zeta ,\bar\zeta ) = L^0_{\rho (z)}(\zeta ,\bar\zeta ) + L^1_{\rho (z)}(\zeta ,\bar\zeta ),
\end{equation}
where
\begin{equation}
\eqqno(levif-2)
L^0_{\rho (z)}(\zeta ,\bar\zeta ) = \frac{1}{4}H^{\rr}_{\rho' (z')}\left( d\Phi_{z}(\zeta ),
d\Phi_{z}(\zeta ) \right)  + \frac{1}{4}H^{\rr}_{\rho' (z')} \left( d\Phi_{z}(i\zeta ),d\Phi_{z}(i\zeta )\right)
\end{equation}
and
\begin{equation}
\eqqno(levif-3)
L^1_{\rho (z)}(\zeta ,\bar\zeta ) = 2\left\langle \nabla \rho' (z'),\d\bar\d \Phi_{z}
(\zeta ,\bar\zeta )\right\rangle = 2Re\left(\dbar\rho' (z'),\d\dbar \Phi_{z}(\zeta ,\bar\zeta)\right) .
\end{equation}
\end{lem}
\proof Here we denote by $d\Phi_{z}$ is the differential of the inverse map $\Phi\deff F^{-1}$ at point $z$, $\nabla\rho' (z')$ the
real gradient of $\rho'$ at $z'$, $\left\langle \cdot,\cdot\right\rangle = Re\left(\cdot ,\cdot \right) $ - the standard 
Euclidean scalar product in $\cc^n$.

Denote by $\nu$ the vector with components $\nu_j=\frac{\d z'_j}{\d z_{ \alpha}}\zeta_{\alpha}$
and by $\mu$ with $\mu_j=\frac{\d z'_j}{\d\bar z_{\alpha}}\bar\zeta_{\alpha}$, \ie 
\[
\nu = \d \Phi_{z}(\zeta ) \text{ } \text{ and } \text{ } \mu =  \dbar \Phi_{z}(\zeta ).
\]
Remark that
\begin{equation}
\nu + \mu = d\Phi_{z}(\zeta ) \text{ } \text{ and } \text{ } i(\nu - \mu) = d\Phi_{z}(i\zeta ).
\eqqno(nu-mu)
\end{equation}
Write

\[
L_{\rho (z)}(\zeta ,\bar\zeta ) = \frac{\d^2\rho }{\d z_{\alpha}\d\bar z_{\beta}}\zeta_{\alpha}\bar\zeta_{\beta} = 
\frac{\d}{\d z_{\alpha}}\left( \frac{\d \rho'}{\d z'_i}\frac{\d z'_i}{\d\bar z_{\beta}} + \frac{\d \rho'}{\d\bar z'_i}
\frac{\d\bar z'_i}{\d\bar z_{\beta}}\right)\zeta_{\alpha}\bar\zeta_{\beta} =
\]
\[
= \frac{\d^2\rho'}{\d z'_i\d z'_j}\frac{\d z'_i}{\d\bar z_{\beta}}\frac{\d z'_j}{\d z_{\alpha}}\zeta_{\alpha}\bar\zeta_{\beta} + \frac{\d^2\rho'}{\d\bar z'_i\d \bar z'_j}\frac{\d\bar z'_i}{\d\bar z_{\beta}}\frac{\d\bar z'_j}{\d z_{\alpha}} \zeta_{\alpha}\bar\zeta_{\beta} + \frac{\d^2 \rho'}{\d z'_i \d\bar z'_j}\left(
\frac{\d z'_i}{\d \bar z_{\beta}}\frac{\d \bar z'_j}{\d z_{\alpha}} + \frac{\d z'_i}{\d z_{\alpha}}\frac{\d\bar z'_j}{\d \bar z_{\beta}}\right) \zeta_{\alpha}\bar\zeta_{\beta} +
\]
\[
 + \left( \frac{\d \rho' }{\d z'_i}\frac{\d^2z'_i}{\d z_{\alpha}\d\bar z_{\beta}} +
 \frac{\d \rho' }{\d\bar z'_i}\frac{\d^2\bar z'_i}{\d z_{\alpha}\d\bar z_{\beta}}\right)
\zeta_{\alpha}\bar\zeta_{\beta} = \frac{\d^2 \rho'}{\d z'_i \d z'_j}\mu_i\nu_j +
\frac{\d^2 \rho'}{\d\bar z'_i \d\bar z'_j}\bar \nu_i\bar \mu_j + \frac{\d^2 \rho'}{\d z'_i \d\bar z'_j}\left[
\mu_i\bar \mu_j+\nu_i\bar \nu_j\right] +
\]
\[
+ \left( \frac{\d \rho '}{\d z'_i}\frac{\d^2z'_i}{\d z_{\alpha}\d\bar z_{\beta}} +
 \frac{\d \rho' }{\d\bar z'_i}\frac{\d^2\bar z'_i}{\d z_{\alpha}\d\bar z_{\beta}}\right)
\zeta_{\alpha}\bar\zeta_{\beta} =  L^0_{\rho (z)}(\nu,\mu) + L^1_{\rho (z)}(\zeta,\bar\zeta)
\]
with
\begin{equation}
\eqqno(2.7)
L^0_{\rho (z)}(\nu , \mu ) = \frac{\d^2 \rho'}{\d z'_i \d z'_j}\nu_i\mu_j +
\frac{\d^2 \rho'}{\d\bar z'_i \d\bar z'_j}\bar \nu_i\bar \mu_j + \frac{\d^2 \rho'}{\d z'_i \d\bar z'_j}\left[
\mu_i\bar \mu_j+\nu_i\bar \nu_j\right]
\end{equation}
and
\begin{equation}
\eqqno(2.8)
L^1_{\rho (z)}(\zeta ,\bar\zeta ) = \left( \frac{\d \rho' }{\d z'_i}\frac{\d^2z'_i}{\d z_{\alpha}\d\bar z_{\beta}} +
 \frac{\d \rho' }{\d\bar z'_i}\frac{\d^2\bar z'_i}{\d z_{\alpha}\d\bar z_{\beta}}\right)
\zeta_{\alpha}\bar\zeta_{\beta}.
\end{equation}

We need to get more information about the structure of both terms $L^0_{\rho}$ and $L^1_{\rho}$ of
the Levi form. Let's prove that the following relation holds
\begin{equation}
L^0_{\rho (z)}(\nu,\mu) = \frac{1}{4}H^{\rr}_{\rho' (z')}\left( \nu +\mu,\nu +\mu\right)  +
\frac{1}{4}H^{\rr}_{\rho' (z')}\left(i(\nu -\mu),i(\nu -\mu)\right).
\eqqno(2.9)
\end{equation}
To see this we make the following change in \eqqref(2.9):

\[
\mu_j = V_j + iW_j, \text{ }
\nu_j = V_j - iW_j.
\]
or
\begin{equation}
V = \frac{1}{2}(\nu+\mu)=\frac{1}{2}d\Phi_{z}(\zeta ), \text{ } 
W = \frac{i}{2}(\nu-\mu)=\frac{1}{2}d\Phi_{z}(i\zeta ).
\eqqno(2.10)
\end{equation}

Then

\[
L^0_{\rho (z)}(\nu ,\mu ) = \frac{\d^2\rho' }{\d z'_i\d z'_j}\left( V_i-iW_i\right) \left( V_j+iW_j\right)
+ \frac{\d^2\rho' }{\d\bar z'_i\d\bar z'_j}\left(\overline{V_i}+i\overline{W_i}\right) \left(\overline{V_j} 
-i\overline{W_j}\right)
+
\]
\[
+ \frac{\d^2\rho' }{\d z'_i\d\bar z'_j}\left[ \left( V_i+iW_i\right)\left(\overline{V_j}-i\overline{W_j}\right)+
\left( V_i-iW_i\right)\left(\overline{V_j}+i\overline{W_j}\right)\right] =
\]
\[
= \frac{\d^2\rho' }{\d z'_i\d z'_j}\left(V_iV_j+W_iW_j\right)  + i\frac{\d^2\rho' }{\d z'_i\d z'_j}\left(V_iW_j- W_iV_j\right) 
+ \frac{\d^2\rho' }{\d\bar z'_i\d\bar z'_j}\left(\overline{V_i}\overline{V_j}+
\overline{W_i}\overline{W_j}\right) +
\]
\[
 + i\frac{\d^2\rho '}{\d\bar z'_i\d\bar z'_j}\left(\overline{W_i}\overline{V_j}-\overline{V_i}\overline{W_j}\right)
+ 2\frac{\d^2\rho' }{\d z'_i\d \bar z'_j}\left( V_i\overline{V_j} + W_i\overline{W_j}\right)
 = H^{\rr}_{\rho' (z')}(V,V) + H^{\rr}_{\rho' (z')}(W,W).
\]
We used the obvious relations $\frac{\d^2\rho' }{\d z'_i\d z'_j}\left(V_iW_j- W_iV_j\right) = 0 = \frac{\d^2\rho' }{\d\bar z'_i\d\bar z'_j}\left(\overline{W_i}\overline{V_j}-\overline{V_i}\overline{W_j}\right)$ and the complex
expression of the real Hessian \eqqref(rc-hess). Therefore 

\begin{equation}
L^0_{\rho (z)}(\nu,\mu) = H^{\rr}_{\rho' (z')}(V,V) + H^{\rr}_{\rho' (z')}(W,W).
\eqqno(2.11)
\end{equation}
From \eqqref(2.10) and \eqqref(2.11) we get the formula \eqqref(levif-2) of the Lemma.

\begin{rema}\rm If the real Hessian of $\rho'$ at $z'$ is positive (resp. non-negative) definite then the component $L^0_{\rho (z)}(\nu ,\mu)$ of the Levi form of $\rho $ at $z=F(z')$ is also positive (resp. non-negative) definite for any $\calc^2$-germ of a
diffeomorphism $F$. 

\end{rema}

Now we turn to $L^1_{\rho }$. Note that in complex notations $\nabla\rho =
\dbar\rho$ as well as that standard Euclidean scalar product $\left\langle \cdot,\cdot\right\rangle $ in $\cc^n$ is equal to the real part of the Hermitian one $\left( \cdot,\cdot \right)$. Therefore from \eqqref(2.8) we get
\[
L^1_{\rho (z)}(\zeta ,\bar\zeta ) = \left( \frac{\d \rho' }{\d z'_i}\frac{\d^2z'_i}{\d z_{\alpha}\d\bar z_{\beta}} + \frac{\d \rho' }{\d\bar z'_i}\frac{\d^2\bar z'_i}{\d z_{\alpha}\d\bar z_{\beta}}\right)
\zeta_{\alpha}\bar\zeta_{\beta} =
\]
\[
= \overline{\left( \dbar\rho' ,\frac{\d^2 z'}{\d\bar z_{\alpha}\d z_{\beta}}
\bar\zeta_{\alpha}\zeta_{\beta}\right)} + \left( \dbar\rho' ,\frac{\d^2 z'}{\d\bar z_{\alpha}\d z_{\beta}}
\bar\zeta_{\alpha}\zeta_{\beta}\right) = \overline{\left( \dbar\rho',\d\dbar \Phi_{z}(\zeta ,\bar\zeta )\right) }
\]
\[
 + \left( \dbar\rho',\d\dbar \Phi_{z}(\zeta ,\bar\zeta )\right) = 2Re\left( \dbar\rho',\d\dbar 
\Phi_{z}(\zeta ,\bar\zeta )\right) =  2\left\langle  \nabla\rho',\d\dbar \Phi_{z}(\zeta ,\bar\zeta )\right\rangle ,
\]
which proves \eqqref(levif-3).

\medskip\qed

\section{Proof of the Theorem}

\noindent{\slsf (\sli $\Longleftrightarrow$ (\sliii }

We start with the proof of the geometric characterization of convex to pseudoconvex mappings given in (\sliii of
the Theorem. By a complex (real) line in $\cc^n$ we mean an $1$-dimensional complex (real) subspace of $\cc^n$. The same for
complex (real) plain. 
Take a complex line $l = \texttt{span} \left\lbrace \zeta \right\rbrace  $ in $T_{z}\cc^n$ and let $\Pi'\subset T_{z'}\cc^n$ 
be the real plain -
image of $l$ under $d\Phi_{z}$, \ie $\Pi' = \texttt{span} \left\lbrace d\Phi_{z}(\zeta), d\Phi_{z}(i\zeta)\right\rbrace $. Let $l':=\d\dbar \Phi_z (l)$ denotes the {\slsf real} (!) line - image of $l$ under the mapping 
\[
\d\dbar \Phi_z : \cc^n_z \to \cc^n_{z'},
\]
defined as 
\[
\zeta \mapsto \d\dbar \Phi (\zeta,\bar\zeta)\deff \frac{\d^2\Phi (z)}{\d z_{\alpha}\dbar z_{\beta}}
\zeta_{\alpha}\bar\zeta_{\beta}.
\]
We consider $l'$ as a real line in $T_{z'}\cc^n$.

\begin{lem}
Suppose that given a diffeomorphism $F: U' \to U$. Then $F$ sends convex quadrics to pseudoconvex hypersurfaces if and only if 
for every $z\in U$ and for all $l':=\d\dbar \Phi_z (l)$ and $\Pi' = d\Phi_{z}(l)$ as above one has $l'\subset \Pi'$.
\end{lem}
\proof  Let us prove the ``only if '' assertion first. We may suppose that $z'=0'$. Take any strictly convex $M'=\{ \rho' (z')=0\} $ defined by a $\calc^2$-function $\rho'$ with positive defined Hessian such that $T_{0'}M'\supset \Pi'$. By $M$ denote the 
image $F(M')$.

Consider the following family of hypersurfaces in $U'$: $M'_t = \{ z: \rho'_t(z') := \rho' (z') + t\left\langle \nabla\rho' (0'),z'\right\rangle = 0\}$, $t\in\rr$, $t\not= -1$. All $M'_t$ are strictly convex (they have the same quadratic part as $M'$), all path through zero and $M'_0=M'$. In addition all $M'_t$ are smooth at zero with $T_{'0}M'_t=T_{0'}M'$ for all $t\not= -1$, because $\nabla\rho'_t(0') = (1+t)\nabla\rho' (0')$.
Moreover, if we take some  $\zeta\in T_{0}^cM$ then $\zeta$ will stay to be complex tangent to all $M_t:=F(M'_t)$ at zero because $T_{0}M_t=dF_{0'}(T_{0'}M'_t)$ is the same for all $t$. From Lemma 2.1 we see that

\begin{equation}
L_{\rho_t(0)}(\zeta ,\bar\zeta ) = L^0_{\rho (0)}(\zeta ,\bar\zeta ) + 2(1+t)\left\langle \nabla\rho' (0'),
\d\dbar \Phi (0) (\zeta ,\bar\zeta )\right\rangle ,
\eqqno(12)
\end{equation}
because $L^0_{\rho_t}(0)(\zeta ,\bar\zeta )=L^0_{\rho }(0)(\zeta,\bar\zeta )$ for all $t$ due to the
fact that coefficients of $L^0_{\rho}$ depend only on the second derivatives of $\rho'$ at $0'$ and on $d\Phi_0$.

\smallskip Suppose $\left\langle \nabla\rho' (0'),\d\dbar \Phi (0) (\zeta ,\bar\zeta )\right\rangle \not= 0$.
Then taking an appropriate $t_0$ we can make $L_{\rho_t(0)}(\zeta ,\bar\zeta ) = 0$ because
$L^0_{\rho (0)}(\zeta ,\bar\zeta )$ do not depend on $t$. Remark that $t_0\not=-1$ because $L^0_{\rho (0)}(\zeta ,\bar\zeta )>0$. 
Now we can deform $M'_t$ letting $t$ run over a neighborhood of $t_0$. $M'_t$
stays strictly convex while the Levi form of $M_t$ changes its sign on the vector $\zeta $. Contradiction 
with assumed property of $F$. Therefore

\begin{equation}
 \left\langle \nabla\rho' (0'),
\d\dbar \Phi (0) (\zeta ,\bar\zeta )\right\rangle = 0,
\eqqno(2.13)
\end{equation}
for every strictly convex $M'=\{ z':\rho' (z')=0\}$ such that $T_{0'}M'\supset \Pi'$. For any vector $v\in T_{0'}\cc^n$ orthogonal to $\Pi'$ we can take a strictly convex hypersurface $M'=\{ z':\rho' (z')=0\}$ such that $\nabla\rho' (0')=v$. Therefore $\d\dbar \Phi (0) (\zeta,\bar\zeta)$ is orthogonal to every such $v$. So $l'=\d\dbar \Phi (0) (l)\subset \Pi'$ and the ``only if '' assertion 
of the lemma is proved.

\medskip To prove the opposite direction take a convex quadric $M'=\{ z': \rho'(z')=0\}$ and set $M=F(M')$. Let $\zeta\in
T_z^cM$. Use again Lemma 2.1. The term $L^0_{\rho (z)}(\zeta , \bar\zeta)$ is clearly positive. The term 
$L^1_{\rho (z)}(\zeta , \bar\zeta)$ is zero because $\d\dbar \Phi_z(\zeta ,\bar\zeta)\in d\Phi_z(\left\langle \zeta \right\rangle)
\subset T_{z'}M'$.

\medskip\qed

\smallskip
Let us reformulate the result obtained as follows (and remark that the equivalence of (\sli and (\sliii in Theorem is proved):

\begin{corol}
If $F$ sends convex quadrics to pseudoconvex hypersurfaces if and only if for every $z\in U$ and every vector 
$\zeta \in T_{z}\cc^n$ the following holds:
\begin{equation}
\d\dbar \Phi_{z} (\zeta ,\bar\zeta ) \in \texttt{span} \left\lbrace d\Phi_{z}(\zeta), d\Phi_{z}(i\zeta)\right\rbrace .
\eqqno(3.3)
\end{equation}

\end{corol}

\smallskip
For the convenience of future references let us formulate the abovementioned statement about holomorphic mappings:
\begin{corol}
A $\calc^2$-difeomorphism $F:U'\to U$ sends pseudoconvex quadrics to pseudoconvex hypersurfaces 
if and only if $F$ is either holomorphic or antiholomorphic.
\end{corol}
\proof  This is well known but still let us give a proof. Suppose, for example, that $\Phi$ is antiholomorphic, 
then $\nu =0$ as defined in \eqqref(nu-mu). Therefore \eqqref(levif-3) tells us that $L^1_{\rho (z)}(\zeta',\bar\zeta')\equiv 0$ in the representation \eqqref(levif-1). Now \eqqref(2.7) shows that
 and   gave us

\[
L_{\rho (z)}(\zeta ,\bar\zeta ) = L_{\rho' (z')}\left( \dbar \Phi_{z}(\zeta ),\overline{\dbar \Phi_{z}(\zeta )}\right) 
\]
for every complex tangent $\zeta$. Conclusion follows.

\smallskip Suppose that, vice versa, $F$ sends pseudoconvex quadrics to pseudoconvex hypersurfaces. 
\eqqref(3.3) shows that $\d\dbar \Phi_{z}(\zeta ,\bar\zeta )$ belongs to the plain
$\texttt{span} \left\lbrace d\Phi_{z}(\zeta), d\Phi_{z}(i\zeta)\right\rbrace $ for all $\zeta\in \cc^n_{z}$. 
And therefore for every $\zeta $ complex tangent to $M=\{ \rho (z)=0\}$ the vector $\d\dbar
\Phi_z(\zeta ,\bar\zeta )$ is tangent to $M'=\{ \rho' (z')=0\}$.
Consequently $L_{\rho (z)}^1(\zeta ,\bar\zeta )\equiv 0$ for any $\rho$. Apply
\eqqref(2.9) to the quadric

\begin{equation}
\rho' (z') = \sum_{j=1}^n(z_j^2+\bar z_j^2 + \eps |z_j|^2) + L(z) + \overline{L(z)}
\eqqno(3.4)
\end{equation}
(where $\eps >0$ and $L$ is a $\cc$-linear form) and get 
\begin{equation}
L^0_{\rho (z)}\left(\zeta ,\bar\zeta \right) = \sum_{j=1}^n(\nu_j\mu_j + \bar\nu_j\bar\mu_j +
\eps |\nu_j|^2 + \eps |\mu_j|^2) = 2Re\left( \sum_{j=1}^n\nu_j\mu_j\right) + \eps (\norm{\nu}^2 + \norm{\mu}^2).
\eqqno(3.5)
\end{equation}
Taking different linear forms $L$ in \eqqref(3.4) we can deploy any $\zeta\in \cc^n$ as a complex tangent 
and therefore, if $\Phi$ is neither holomorphic no antiholomorphic, then we see from \eqqref(nu-mu) that 
 $\nu$ and $\mu$ can be taken arbitrary. But for arbitrary taken $\nu$ and $\mu$ \eqqref(3.5) cannot
be positive. Contradiction.

\smallskip\qed

\smallskip\noindent{\slsf (\sliii $\Longleftrightarrow$ (\slii }

\smallskip We shall need the following linear algebra lemma. Let $V$ and $W$ be 
$\cc$-linear spaces. We suppose that on $V$ some Hermitian scalar  product
$\left(  \cdot , \cdot \right)$ is fixed. Let $B(\zeta,\bar\eta): V\times V\to W$ be a sesquilinear map.
Its trace is defined as $TrB=\sum_{\alpha}B(e_{\alpha},\bar e_{\alpha})$ for an orthonormal frame in 
$\left( V, \left( \cdot , \cdot\right) \right) $. Let, furthermore $C:V\to W$ be an $\rr$-linear isomorphism.
Denote by $C^{1,0}$ (resp. $C^{0,1}$) the complex linear (resp. antilinear) part of $C$.

\begin{lem}
The following properties ob the pair $(B,C)$ are equivalent:

\begin{equation}
B(\zeta,\bar\zeta) \in \texttt{span}\{ C(\zeta), C(i\zeta)\} \texttt{ } \text{ for all } \texttt{ } \zeta\in V.
\eqqno(prop1)
\end{equation}

\begin{equation}
\eqqno(prop2)
B(\zeta ,\bar\eta)=\left(C^{-1}(TrB),\eta\right)C^{1,0}(\zeta) +  \left(\zeta,C^{-1}(TrB)\right)C^{0,1}(\eta)
\texttt{ } \text{ for all } \texttt{ } \zeta , \eta \in V.
\end{equation}

\end{lem}
\proof Define the induced quadratic map $A:V\to V$ as $A(\zeta,\bar\zeta)=C^{-1}\circ B (\zeta,\bar\zeta)$.
Note that $A$ is not sesquilinear in general. Note that the image of every complex 
line in $V$ under a quadratic map is a real line. 

Write \eqqref(prop1) in the form $A(\zeta,\bar\zeta)=k(\zeta)\cdot\zeta$, where $k$ is a complex
valued function. One readily sees that $k(\lambda\zeta)=\bar\lambda
k(\zeta)$. The polarization equality for $A$
\[
A(\zeta + \eta,\bar\zeta + \bar\eta) + A(\zeta - \eta,\bar\zeta - \bar\eta) = 2A(\zeta,\bar\zeta) +
2A(\eta,\bar\eta)
\]
gives
\[
k(\zeta +\eta)(\zeta +\eta) + k(\zeta -\eta)(\zeta -\eta) =
2k(\zeta)\zeta + 2k(\eta)\eta ,
\]
or, for complex independent
vectors
\[
k(\zeta +\eta)+k(\zeta -\eta) = 2k(\zeta) \text{ }\text{and}\text{ }k(\zeta +\eta)-k(\zeta -\eta) = 2k(\eta),
\]
which implies additivity of $k$: $k(\zeta +\eta)=k(\zeta)+k(\eta)$
for complex independent $\zeta,\eta$ and, therefore for all. So $k$
is an antilinear form on $V$ and by Ries representation we obtain a
vector  $v$ such that $k(\zeta) = (v,\zeta)$ for all $\zeta \in V$
and therefore $A(\zeta,\bar\zeta)=(v,\zeta)\zeta $ and consequently

\begin{equation}
B(\zeta,\bar\zeta)=C\left( (v,\zeta)\zeta \right) \texttt{ } \text{ for all } \texttt{ } \zeta\in V.
\end{equation}
Furthermore,
\[
B(\zeta,\bar\eta)+B(\eta,\bar\zeta) = B(\zeta +\eta,\bar\zeta
+\bar\eta) - B(\zeta,\bar\zeta) - B(\eta,\bar\eta) = C\left( (v,\zeta
+\eta)(\zeta +\eta)\right)  -
\]
\[
- C\left( (v,\zeta)\zeta\right)  - C\left( (v,\eta)\eta\right)  = C\left( (v,\eta)\zeta\right)  + 
C\left( (v,\zeta)\eta \right) .
\]
and 
\[
-iB(\zeta,\bar\eta)+iB(\eta,\bar\zeta) = B(\zeta +i\eta,\bar\zeta
-i\bar\eta) - B(\zeta,\bar\zeta) - B(\eta,\bar\eta) = C\left( (v,\zeta
+i\eta)(\zeta +i\eta)\right)  -
\]
\[
- C\left( (v,\zeta)\zeta\right)  - C\left( (v,\eta)\eta\right)  = C\left(-i(v,\eta)\zeta\right)  + 
C\left(i(v,\zeta)\eta \right) .
\]
Therefore 
\[
2B(\zeta,\bar\eta) = C\left( (v,\eta)\zeta\right) -iC\left( v,\eta)\zeta\right) + C\left( (v,\zeta)\eta\right) 
+iC\left( i((v,\zeta)\eta\right).
\]
So we obtain
\begin{equation}
B(\zeta,\bar\eta) = (v,\eta)C^{1,0}(\zeta )+ (\zeta,v)C^{0,1}(\eta ).
\eqqno(20)
\end{equation}
Set in \eqqref(20) $\zeta=\eta= e_{\alpha}$. Then 
\[
TrB = \sum_{\alpha}B(e_{\alpha},\bar e_{\alpha}) = C^{1,0}\left( \sum_{\alpha}(v,e_{\alpha})\right) +
 C^{0,1}\left( \sum_{\alpha}(v,e_{\alpha})\right) = C(v).
\]
Therefore $v=C^{-1}\left( TrB\right) $ and \eqqref(prop2) is established.

\medskip The opposite implication is easy, because \eqqref(prop2) tells, if $\eta$ is taken to be equal to
$\zeta$,  that 
\[
B(\zeta,\bar\zeta) = aC^{1,0}(\zeta) + \bar a C^{0,1}(\zeta) = a\frac{1}{2}\left( C(\zeta) - i C(i\zeta)\right)
+ \bar a\frac{1}{2}\left( C(\zeta) + i C(i\zeta)\right) = 
\]
\[
= \re a \cdot C(\zeta) + \im a \cdot C(i\zeta) \in \texttt{span}\{ C(\zeta), C(i\zeta)\}.
\]

\medskip\qed

We apply this lemma for $B=\d\dbar \Phi_{z}:T_{z}\cc^n\to T_{z'}\cc^n$, $C=dF_{z'}^{-1}:T_{z'}\cc^n\to T_{z}\cc^n$ 
and, as a result $A=dF_{z'}\circ \d\dbar \Phi_{z}:T_{z}\cc^n\to T_{z}\cc^n$ for
every $z=F(z')$ and get

\begin{corol}
A $\calc^2$-diffeomorphism $F$ sends convex quadrics to pseudoconvex hypersurfaces
if and only if 
\begin{equation}
\d\dbar \Phi_{z}(\zeta ,\bar\eta ) = \left( dF_{z'}\left( Tr\d\dbar \Phi_{z}\right) ,\eta\right)
\d \Phi(\zeta) + \left( \zeta, dF_{z'}\left( Tr\d\dbar \Phi_{z}\right)\right)
\dbar \Phi(\eta)
\end{equation}
for all $z=F(z')$ and all $\zeta ,\eta \in T_{z}\cc^n$.
\end{corol}

And this is equivalent to \eqqref(syst1). Theorem is proved.

\ifx\undefined\bysame
\newcommand{\bysame}{\leavevmode\hbox to3em{\hrulefill}\,}
\fi

\def\entry#1#2#3#4\par{\bibitem[#1]{#1}
{\textsc{#2 }}{\sl{#3} }#4\par\vskip2pt}

\end{document}